\documentclass{article}
\usepackage{amsmath}

\input{tcilatex}

\begin{document}

\begin{center}
\bigskip 

\textbf{OPTIMAL CONTROL OF DAMS USING }$\mathbf{P}_{\lambda ,\tau }^{M}$%
\textbf{\ POLICIES AND PENALTY COST WHEN THE INPUT PROCESS IS }

\textbf{AN INVERSE GAUSSIAN PROCESS}

\ 

\bigskip

Mohamed Abdel-Hameed

Department of Statistics

College of Business and Economics

United Arab Emirates University
\end{center}

\bigskip

\begin{center}
\textbf{ABSTRACT }
\end{center}

We consider $P_{\lambda ,\tau }^{M}\ $ policy of a dam in which the water
input is an inverse Gaussian process. The release rate of the water is
changed from $0$ to $M$ and from $M$ to $0$ $(M>0)$ at the moments when the
water level up crosses levels $\lambda \ $ and down crosses level $\tau $ $%
(\tau <\lambda )$, respectively. We determine the resolvent of the dam
content and compute the total discounted as well as the long-run-average
cost. We also find the stationary distribution of the dam content.

\bigskip

\smallskip

\begin{center}
I. \textbf{INTRODUCTION}
\end{center}

\noindent Adel-Hameed (2000) discuss the optimal control of a dam using P$%
_{\lambda ,\tau }^{M}$\ policies, using total discounted cost as well as the
long -run average cost. He assumes that the water input is a compound
Poisson process with a positive drift. The release rate is zero until the
water reaches level $\lambda ,\;$\ then it is released at rate $M$ until it
reaches level $(\tau <\lambda )$, once the water\ reaches level $\tau \;$\
the release rate remains zero until level $\lambda \;$\ is reached again and
the cycle is repeated. At any time , the release rate can be increased from
zero to M with a starting cost$\;K_{1}M\;$, or decreased from$\;M$\ to zero
with a closing cost $K_{2}M\;.$\ Moreover, for each unit of output, a reward 
$R\;$is received.\ Furthermore, there is a penalty cost which accrues at a
rate $g,$\ where $g$ is a bounded measurable function on the state space.
When the release rate is zero, the water content is denoted by \ $%
I\;=(I_{t}).\;$Let$\;W_{\lambda }\;$be the first passage time of the process 
$I$\ through the boundary $\lambda .$\ When the release rate is M the
content process is denoted by $I^{\ast }=(I_{t}^{\ast }),$\ the process $%
I^{\ast }\ $is a strong Markov process.\ Let$\;W_{\tau }^{\ast }\;$be the
first passage of process $I^{\ast }$\ to the boundary $\tau $. Over time,
the content process is obtained by hitching independent copies of the
processes $I$\ and $I^{\ast }\;$together. The content process is best
described by the bivariate process $B=(Z,R)$, where $Z=(Z_{t}),R=(R_{t})\;$%
describe the dam content, and the release rate respectively. We define the
following sequence of stopping times : 
\begin{eqnarray*}
\overset{\symbol{94}}{T}_{0} &=&\inf \{t\geq 0:Z_{t}\geq \lambda
\},\;\;\;\;\;\overset{\ast }{T}_{0}=\inf \{t\geq \overset{\symbol{94}}{T}%
_{0}:Z_{t}=\tau \} \\
\overset{\symbol{94}}{T}_{n} &=&\inf \{t\geq \overset{\ast }{T}%
_{n-1}:Z_{t}\geq \lambda \},\;\overset{\ast }{T_{n}}=\inf \{t\geq \overset{%
\symbol{94}}{T}_{n}:Z_{t}=\tau \},\;for\;n\geq 1.
\end{eqnarray*}%
It follows that the process$\;B$\ is a delayed regenerative process with the
regeneration points being the $\overset{\ast }{T}_{n},\;n=0,1,...$\ $.$\ $\;$%
The penalty cost rate function is defined as follows 
\[
g(z,r)=\left\{ 
\begin{array}{ll}
g(z) & \;(z,r)\in (0,\lambda )\times \{0\} \\ 
g^{\ast }(z) & \;(z,r)\in (\tau ,\infty )\times \{M\}%
\end{array}%
\right. 
\]%
where $g:(0,\lambda ):\rightarrow R_{+},$\ $g^{\ast }:(\tau ,\infty
):\rightarrow R_{+}\;$are bounded measurable functions.

For any process $Y=(Y_{t},t\geq 0\},\;$and functional f, $E_{y}(f)\;$denotes
the expectation of\ $f$\ conditional on $Y_{0}=y,\;$and\ $P_{y}(A)\;$denotes
the corresponding probability measure. Throughout, we let $R_{+}=[0,\infty
),\;$\ $N=\{1,2,...\}$, $N_{+}=\{0,1,...\}$, and $I_{A}(\;)\;$be the
indicator function of any set $A$. Let $C_{g}^{\alpha }(0,x,\lambda ),\;$be
the expected discounted penalty costs during the interval $(0,W_{\lambda
})\; $starting at x$,\;$and $C_{g^{\ast }}^{\alpha }(M,y,\tau )\;$be the
expected discounted penalty costs , starting at y, during the interval ($%
0,W_{\tau }^{\ast })$\ .\ Furthermore, let $C_{g}(0,x,\lambda ),C_{g^{\ast
}}(M,y,\tau )\;$be the expected non-discounted penalty costs during the same
intervals. respectively.\ It follows that 
\begin{eqnarray*}
C_{g}^{\alpha }(0,x,\lambda ) &=&E_{x}\int_{0}^{W_{\lambda }}e^{-\alpha
t}g(I_{t})dt,\;C_{g^{\ast }}^{\alpha }(M,y,\tau )\;=E_{y}\int_{0}^{W_{\tau
}^{\ast }}e^{-\alpha t}g(I_{t}^{\ast })dt, \\
C_{g}(0,x,\lambda ) &=&E_{\tau }\int_{0}^{W_{\lambda
}}g(I_{t})dt,\;\;C_{g^{\ast }}(M,y,\tau )\;=E_{y}\int_{0}^{W_{\tau }^{\ast
}}g^{\ast }(I_{t}^{\ast })dt
\end{eqnarray*}
Bae \textit{et al.} (2003) consider the average cost case of the above
model, when the dam has a finite capacity. Abdel-Hameed and Nakhi (2006)
treat the case where the water input is a diffusion process. The release
rate depends on the water content. Dohi \textit{et al.} (1995) consider the
case where the water must be released at a fixed time $T$. They assume that
the input is a wiener process. In this paper, we discuss the case where the
water input is an inverse Gaussian process. The assumption that the water
input is of that type is more realistic than the cases where it is assumed
that it is a Wiener process. This is true because the inverse Gaussian
process has increasing sample paths. We determine the total discounted as
well as the long -run average costs. In Section\ II,\ we discuss the
resolvent operators of the processes of interest. In Section III we obtain
formulas for the cost functional using the total discounted as well as the
long-run average cost criteria.

\begin{center}
\textbf{\ }$\;$\bigskip

II. \textbf{THE DAM CONTENT PROCESS AND ITS CHARACTERISTICS}
\end{center}

\noindent Assume that the water input in the dam $I=$ $\{I_{t},t\geq 0\}$ is
an inverse Gaussian process, with parameters $\mu ,\ \sigma ^{2}\ \geq 0\ $.
We let $p(t,x,y\ )\ $be the probability transition function of the process $%
I\ $given $I_{0}=x$, and $p(t,y\ )\ $be its probability transition function
given $I_{0}=0$. Since the process $I\ $is additive, we have that for $y\geq
0\ $and$\ x\ \geq 0\ $ 
\begin{eqnarray*}
p(t,x,y\ )\ \ &=&\ p(t,y-x\ ) \\
&=&\left\{ 
\begin{array}{l}
\frac{t}{\sigma \sqrt{2\pi (y-x)^{3}}}\exp \{-\frac{[\mu (y-x)-t]^{2}}{%
2(y-x)\sigma ^{2}}\},\ y>x \\ 
0\ \ \ \ \ \ \ \ \ \ \ \ \ \ \ \ \ \ \ \ \ \ \ \ \ \ \ \ \ \ \ \ \ \ \ \ \ \
\ \ \ \ \ ,y\leq x\text{.}%
\end{array}%
(1)\right.
\end{eqnarray*}%
It is known that $I$ is a pure-jump process, with state space $(0,\infty ),$%
and with jump measure $\nu $ concentrated on $(0,\infty )$ given by 
\[
\nu (dy)=\frac{1}{\sigma \sqrt{2\pi }}\frac{\exp (-y\mu ^{2}/2\sigma ^{2})}{%
y^{3/2}} 
\]

It follows that for $\alpha \geq 0$ 
\[
Ee^{-\alpha I_{t}}=e^{-\frac{t}{\sigma ^{2}}\{\sqrt{2\alpha \sigma ^{2}+\mu
^{2}}-\mu \}}\ \ 
\]%
Direct differentiation of the above function gives, $EX_{t}=\frac{t}{\mu }$,
and $Var(X_{t})=\frac{t\sigma ^{2}}{\mu ^{3}}$.

$\;$To evaluate the cost functional and other parameters of the process
during the first part of a cycle, we define the Levy process killed at $%
\lambda ,\;$as follows$:$%
\[
X=\{I_{t},t<W_{\lambda }\} 
\]%
Let $U_{\alpha }\;$be the resolvent operator of the process $I,$\ defined
for every bounded function $f$ \ by $U_{\alpha
}f(x)=\dint\limits_{0}^{\infty }f(x+y)U_{\alpha }(dy)$. \ Then $U_{\alpha }\ 
$is the unique solution of\bigskip\ the equation

\[
\dint\limits_{0}^{\infty }f(x+y)U_{\alpha
}(dy)=E_{x}\dint\limits_{0}^{\infty }e^{-\alpha t}f(X_{t})dt\text{.} 
\]

It follows with $f_{\beta }(x)=exp(-\beta x),\beta \geq 0$, that

\[
U_{\alpha }f_{\beta }(0)=\frac{\sigma ^{2}}{\alpha \sigma ^{2}+\{\sqrt{%
2\beta \sigma ^{2}+\mu ^{2}}-\mu \}.} 
\]

\bigskip\ \ Throughout we let $\varphi _{_{Z}}(.)$ as the standard normal
density function, $\func{erf}()\ $and $\func{erf}c()\ $be the well known
error and complimentary error functions, respectively. Inverting the above
function w.r.t $\beta \ $we have

\bigskip

\ 
\begin{eqnarray*}
U_{\alpha }(dy) &=&\frac{\sigma }{\sqrt{y}}\varphi _{_{Z}}(\sqrt{y}\mu
/\sigma )dy+(\frac{\mu -\alpha \sigma ^{2}}{2})e^{\alpha y(\frac{\alpha
\sigma ^{2}}{2}-\mu )}\func{erf}c\{\sqrt{y}\frac{\alpha \sigma ^{2}-\mu }{%
\sqrt{2\sigma ^{2}}}\}dy\ \ \text{(2).\ \ } \\
&=&u_{\alpha }(y)dy
\end{eqnarray*}

where

\[
u_{\alpha }(y)=\frac{\sigma }{\sqrt{y}}\varphi _{_{Z}}(\sqrt{y}\mu /\sigma
)+(\frac{\mu -\alpha \sigma ^{2}}{2})e^{\alpha y(\frac{\alpha \sigma ^{2}}{2}%
-\mu )}\func{erf}c\{\sqrt{y}\frac{\alpha \sigma ^{2}-\mu }{\sqrt{2\sigma ^{2}%
}}\}. 
\]

It follows that, for $x\leq \lambda ,$ 
\begin{eqnarray*}
C_{g}^{\alpha }(0,x,\lambda ) &=&\int_{0}^{\lambda -x}g(x+y)U_{\alpha }(dy),
\\
C_{g}(0,x,\lambda ) &=&\int_{0}^{\lambda -x}g(x+y)U_{0}(dy).
\end{eqnarray*}%
From Equation (8) of [7], it follows that for $x\leq \lambda ,$

\bigskip

\begin{eqnarray*}
E_{x}(\exp (-\alpha W_{\lambda })) &=&\alpha U_{\alpha }I_{[\lambda ,\infty
)}(x)\text{. } \\
&=&\frac{\alpha \sigma ^{2}-\mu }{\alpha \sigma ^{2}-2\mu }e^{\alpha
(\lambda -x)(\frac{\alpha \sigma ^{2}}{2}-\mu )}\func{erf}c(\sqrt{\lambda -x}%
\frac{\alpha \sigma ^{2}-\mu }{\sqrt{2\sigma ^{2}}}) \\
&&-\frac{\mu }{\alpha \sigma ^{2}-2\mu }\func{erf}c(\frac{\sqrt{\lambda -x}%
\mu }{\sqrt{2\sigma ^{2}}})\ \ \ \ \ \ \ \ \ \ \ \ \ \ \ \ \text{\ (3)}
\end{eqnarray*}%
where the last equation follows by integrating $U_{\alpha }(dy)\ $over the
interval $[\lambda ,\infty )\ $(we omit the proof)$.$

\bigskip

It follows that, for $x\leq \lambda $,\ the distribution function of $%
W_{\lambda }\ $(denoted by $F_{W_{\lambda }}()$) is given by

\bigskip 
\[
F_{W_{\lambda }}(t)=\frac{1}{2}\func{erf}c\{\frac{(\lambda -x)\mu -t}{\sqrt{%
2\sigma ^{2}}}\}-\frac{1}{2}e^{2\mu t/\sigma ^{2}}\func{erf}c\{\frac{%
(\lambda -x)\mu +t}{\sqrt{2\sigma ^{2}}}\},\ t\geq 0\text{.} 
\]%
Furthermore, for $x\leq \lambda $

\begin{eqnarray*}
E_{x}(W_{\lambda }) &=&U_{0}I_{[0,\lambda )}(x)\;\;\; \\
&=&\sigma \int_{0}^{\lambda -x}\frac{1}{\sqrt{y}}\varphi _{_{Z}}(\sqrt{y}%
\frac{\mu }{\sigma })dy+\frac{\mu }{2}\int_{0}^{\lambda -x}\func{erf}c(-%
\sqrt{\frac{y}{2}}\frac{\mu }{\sigma })\ dy. \\
&=&\frac{(\lambda -x)\mu }{2}+\sigma \sqrt{\lambda -x\ }\varphi _{_{Z}}(%
\sqrt{\lambda -x}\frac{\mu }{\sigma })+\frac{(\lambda -x)\mu ^{2}+\sigma ^{2}%
}{2\mu }\func{erf}(\sqrt{\frac{\lambda -x}{2}}\frac{\mu }{\sigma })\text{ \ }%
(4)
\end{eqnarray*}%
where the last equation follows from the equation before last upon tedious
calculations which we omit.

\bigskip To derive $C_{g^{\ast }}^{\alpha }(M,y,\tau ),C_{g^{\ast
}}(M,y,\tau ),E_{y}(\exp (-\alpha W_{\tau }^{\ast })),\;$and\ $E_{y}$($%
W_{\tau }^{\ast })$.$\;\ $We first note that $I\ ^{\ast }\ $is a Levy
process and using Doob's optional sampling theorem\ we have the following
result%
\[
E_{x}[e^{-\alpha W_{\tau }^{\ast }}]=e^{-(x-\tau )\eta (\alpha )}\text{ .\ }%
(5) 
\]%
where $\eta (\alpha )$ is the unique increasing solution of the equation 
\[
M\eta (\alpha )=\alpha +\frac{\sqrt{2\eta (\alpha )\sigma ^{2}+\mu ^{2}}-\mu 
}{\sigma ^{2}}.(6) 
\]

It can be seen that the permissible solution of this equation is (we omit
the proof)

\bigskip

\[
\eta (\alpha )=\frac{\alpha }{M}+\frac{(1-M\mu )+\sqrt{2\alpha M\sigma
^{2}+(1-M\mu )^{2}}}{M^{^{2}}\sigma ^{2}}. 
\]%
It follows that, for any $x\geq \tau ,$

\[
\Pr_{x}\{W_{\tau }^{\ast }\,<\infty \}\text{ =}\left\{ 
\begin{array}{ll}
1 & \ \ \text{if\ }\mu M>1 \\ 
\ e^{\frac{-2(x-\tau )(1-\mu M)}{M^{2}\sigma ^{2}}} & \ \ \text{if }\mu
M\,\,\leq 1\text{.}%
\end{array}%
\right. 
\]%
It is also found that the probability density function of $W_{\tau }^{\ast
}\ \ (f_{W_{\tau }^{\ast }}(.))\ $is equal to zero\ for $t<\frac{(x-\tau )}{M%
}$, and for $t\geq \frac{(x-\tau )}{M},$

\[
f_{W_{\tau }^{\ast }}(t)=\frac{(x-\tau )}{\sigma \sqrt{2\pi (Mt-(x-\tau
))^{3}}}\exp \{-\frac{((M\mu -1)t-\mu (x-\tau ))^{2}}{2(Mt-(x-\tau ))\sigma
^{2}}\}. 
\]%
Furthermore,

\begin{eqnarray*}
E_{x}W_{\tau }^{\ast } &=&\frac{(x-\tau )\mu }{(\mu M-1)}\text{ \ \ if }\mu
M>1,\ \ (7) \\
&=&\infty \;\;\;\;\;\;\;\;\;\;\ \;\text{if }\mu M\leq 1,
\end{eqnarray*}%
and$\;$

\begin{eqnarray*}
Var_{x}(W_{\tau }^{\ast }) &=&\frac{(x-\tau )M^{2}\sigma ^{2}}{(\mu M-1)^{3}}%
\;\;\text{\ if \ }\mu M\ >1, \\
&=&\infty \;\;\;\;\;\;\;\;\;\;\;\;\;\ \ \ \ \text{if \ }\mu M\leq 1\text{.}
\end{eqnarray*}

\ \ \ \ \ \ \ \ \ \ We now define, the killed process 
\[
X^{\ast }=\{\overset{\ast }{I}_{t},t\leq W_{\tau }^{\ast }\} 
\]%
It can be shown that the process $X^{\ast }\;$is a strong Markov process.
Furthermore, it has state space $(\tau ,\infty ).$\ Starting at $\;x$\ $\in $%
\ $(\tau ,\infty )$, let $f(x,y,t)\ $be the transition probability function
of\ $\overset{\ast }{I}$. Let $U_{\alpha }^{\ast }\;$be the resolvent
operator of the process $X^{\ast }$,$\;$ it follows that for $x\in (\tau
,\infty )$,\ 

\[
U_{\alpha }^{\ast }(dy-x)=[p_{\alpha }^{\ast }(y-x)-\exp (-(x-\tau )\eta
(\alpha ))p_{\alpha }^{\ast }(y-\tau )]dy\text{.} 
\]%
where for any $z$ ,

\[
p_{\alpha }^{\ast }(z)=\int_{0}^{\infty }\exp (-\alpha t)p(t,z+Mt)dt\text{,} 
\]%
where $p(t,x)\ $is the transition probability of the process $I$, \ starting
at zero, defined before.

\bigskip

Furthermore, for $x,y\geq \tau $,

$\ $

$\bigskip $%
\[
f(x,y,t)=p(t,y-x+Mt)-\int\limits_{0}^{t}p(t-s,y-\tau +M(t-s))f_{W_{\tau
}^{\ast }}(s)ds\text{.} 
\]%
Thus, for $x\geq \tau ,$%
\[
C_{g^{\ast }}^{\alpha }(M,x,\tau )=\int_{\tau }^{\infty }g^{\ast
}(y)U_{\alpha }^{\ast }(dy-x)\text{,} 
\]%
and 
\[
C_{g^{\ast }}(M,x,\tau )=\int_{\tau }^{\infty }g^{\ast }(y)U_{0}^{\ast
}(dy-x). 
\]

\indent Now we need to compute the joint distribution function of the pair $%
(W_{\lambda },I_{W_{\lambda }})\ $given $I_{0}=0$, denoted by $f_{0}(t,x)$.
We define $v_{\alpha }(x)=u_{\alpha }(x)I\{x>\lambda \}$ and $v_{\alpha
}^{0}(x)=u_{\alpha }(x)I\{x=\lambda \}$. For any function $g\ $we \ let $%
L_{\beta }(g)\ $be its Laplace transform with respect to $\beta .\ $From
Equation (8) P. 2067 \ of [7] we have, for $\alpha \geq 0$, $\beta \geq 0$

\bigskip 
\begin{eqnarray*}
E[e^{-\alpha W_{\lambda }-\beta I_{W_{\lambda }}}] &=&[\alpha +\frac{\sqrt{%
2\beta \sigma ^{2}+\mu ^{2}}-\mu }{\sigma ^{2}}]L_{\beta }(v_{\alpha }). \\
&=&\alpha L_{\beta }(v_{\alpha })+\frac{\sqrt{2\beta \sigma ^{2}+\mu ^{2}}%
-\mu }{\sigma ^{2}}L_{\beta }(v_{\alpha }). \\
&=&\alpha L_{\beta }(v_{\alpha })+\frac{2\beta }{\sqrt{2\beta \sigma
^{2}+\mu ^{2}}-\mu }L_{\beta }(v_{\alpha })-\frac{2\mu }{\sigma ^{2}}%
L_{\beta }(v_{\alpha }) \\
&=&\alpha L_{\beta }(v_{\alpha })+2L_{\beta }(u_{0}\ast v_{\alpha
}^{o})++2L_{\beta }(u_{0})L_{\beta }(v_{\alpha }^{^{\prime }})-\frac{2\mu }{%
\sigma ^{2}}L_{\beta }(v_{\alpha }) \\
&=&\alpha L_{\beta }(v_{\alpha })+2L_{\beta }(u_{0}\ast v_{\alpha
}^{o})++2L_{\beta }(u_{0}\ast v_{\alpha }^{^{\prime }})-\frac{2\mu }{\sigma
^{2}}L_{\beta }(v_{\alpha })\text{. \ }\ (8)
\end{eqnarray*}

Inverting the above function with respect to $\alpha ,\beta ,$ we have

\[
f_{0}(t,x)=\frac{\partial }{\partial t}p(t,x)+2\{p_{t}(\lambda
)u_{0}(x-\lambda )+\int_{\lambda }^{x}u_{0}(x-y)\frac{\partial }{\partial y}%
p(t,y)dy-\frac{\mu }{\sigma ^{2}}p(t,x)\}]I\{t\geq 0,x>\lambda \}\ \ \ (9) 
\]

To find the marginal pdf of $I_{W_{\lambda }}$, we define

\[
L_{\beta }(\lambda )\overset{def}{=}E_{0}(e^{^{--\beta I_{W_{\lambda }}}})%
\text{.} 
\]%
Letting $\alpha \downarrow 0\ \ $in$\ $Equation$\ $(8) we get

\bigskip

\[
L_{\beta }(\lambda )=2[L_{\beta }(u_{0}\ast v_{0}^{0})+L_{\beta }(u_{0}\ast
v_{0}^{\prime })-\frac{\mu }{\sigma ^{2}}L_{\beta }(v_{0})]. 
\]%
It follows that

\bigskip 
\begin{eqnarray*}
E_{x}(I_{W_{\lambda }}) &=&\frac{U_{0}I_{(0,\lambda -x]}(0)}{\mu }\;\ (10) \\
&=&\frac{E_{x}(W_{\lambda })}{\mu }
\end{eqnarray*}%
Furthermore, the pdf of $I_{W_{\lambda }}$, given $I_{0}=0\ $(denoted by $%
f_{I_{W_{\lambda }}}()$) is given by

\[
f_{I_{W_{\lambda }}}(x)=2[u_{0}(\lambda )u_{0}(x-\lambda )+\int_{\lambda
}^{x}u_{0}(x-y)u_{0}^{^{\prime }}(y)dy-\frac{\mu }{\sigma ^{2}}%
u_{0}(x)\}]I\{x>\lambda \}. 
\]

\bigskip

\bigskip

\begin{center}
III. \textbf{\ THE TOTAL DISCOUNTED AND LONG RUN-AVERAGE COSTS AND THE
STATIONARY DISTRIBUTION OF THE DAM CONTENT}
\end{center}

\noindent We now discuss the computations of the cost functional using the
total discounted cost as well as the long-run average cost. Let $W\ $be \ \
the length of the first cycle, i.e., $W$=$\overset{\ast }{T}_{1}-\overset{%
\ast }{T}_{0},\;$and $C_{\alpha }(x)\;$be the expected cost during the
interval $[0,\overset{\ast }{T}_{0}),\;$when Z$_{0}=x.\;$Then,\ it follows
that the total discounted cost associated with an $P_{\lambda ,\tau }^{M}\;$%
policy is given by

\[
C_{\alpha }(\lambda ,\tau )=C_{\alpha }(x)+\frac{E_{x}(\exp (-\alpha \overset%
{\ast }{T}_{0})E_{\tau }C_{\alpha }(1)}{1-E_{\tau }(\exp (-\alpha W))} 
\]%
where\ $C_{\alpha }(1)\ $is the total discounted cost during the interval $%
(0,W).$\ For $x\leq \lambda $,$\ $ we have

\begin{eqnarray*}
C_{\alpha }(x) &=&M\{K_{2}+K_{1}E_{x}(e^{-\alpha W_{\lambda
}})-RE_{_{x}}\int_{W_{\lambda }}^{\overset{\ast }{T}_{0}}e^{-\alpha t}dt\} \\
&&+\int_{0}^{\overset{\ast }{T}_{0}}e^{-\alpha t}g(Z_{t},R_{t})dt \\
&=&M\{K_{2}+K_{1}E_{x}(e^{-\alpha W_{\lambda }})-RE_{_{x}}\int_{W_{\lambda
}}^{\overset{\ast }{T}_{0}}e^{-\alpha t}dt\} \\
&&+E_{_{x}}\int_{0}^{W_{\lambda }}e^{-\alpha
t}g(I_{t})dt+E_{_{x}}\int_{W_{\lambda }}^{\overset{\ast }{T}_{0}}e^{-\alpha
t}g^{\ast }(I_{t}^{\ast })dt \\
&=&M\{K_{2}+K_{1}E_{x}(e^{-\alpha W_{\lambda }})-RE_{_{x}}\{e^{-\alpha
W_{\lambda }}E_{I_{W_{\lambda }}}\int_{0}^{\overset{\ast }{W}_{\tau
}}e^{-\alpha t}dt\} \\
&&+C_{g}^{\alpha }(0,\tau ,\lambda )+E_{x}[e^{-\alpha W_{\lambda
}}E_{I_{W_{\lambda }}}\int_{0}^{\overset{\ast }{W}_{\tau }}e^{-\alpha
t}g^{\ast }(I_{t}^{\ast })dt) \\
&=&M\{K_{2}+K_{1}E_{x}(e^{-\alpha W_{\lambda }})\}+C_{g}^{\alpha }(0,\tau
,\lambda )+E_{x}[e^{-\alpha W_{\lambda }}C_{g^{\ast }-RM}^{\alpha
}(M,I_{W_{\lambda }},\tau )],\ (11)
\end{eqnarray*}%
where the third equation above follows from the second equation upon
conditioning on the sigma algebra generated by $(W_{\lambda },I_{W_{\lambda
}})$.Furthermore, 
\[
E_{\tau }C_{\alpha }(1)=C_{\alpha }(\tau )\text{,} 
\]%
where $C_{\alpha }(\tau )\ $is obtained from Equation (10)\ upon
substituting $\tau \ $for $x$.

Throughout the remainder of this section we let $\varphi (\alpha )=$ $\frac{%
\sqrt{2\alpha \sigma ^{2}+\mu ^{2}}-\mu }{\sigma ^{2}}$.$\ $Now, for $%
x<\lambda $

\begin{eqnarray*}
E_{x}(e^{-\alpha \overset{\ast }{T}_{0}}) &=&E_{x}[e^{-\alpha W_{\lambda
}}E_{I_{W_{\lambda }}}(e^{-\alpha \overset{\ast }{W}_{\tau }})] \\
&=&E_{x}[e^{-\alpha W_{\lambda }}e^{-\eta (\alpha )(I_{W_{\lambda }}-\tau
)\eta (\alpha )}] \\
&=&E_{0}[e^{-\alpha W_{\lambda -x}}e^{-\eta (\alpha )(I_{W_{\lambda
-x}}+x-\tau )}] \\
&=&e^{-\eta (\alpha )(x-\tau )}E_{0}[e^{-\alpha W_{\lambda -x}-\eta (\alpha
)(I_{W_{\lambda -x}})}] \\
&=&(\alpha +\varphi (\eta (\alpha ))e^{-\eta (\alpha )(x-\tau
)}\int_{[\lambda -x,\infty )}e^{-z\eta (\alpha )}U_{\alpha }(dz) \\
&=&M\eta (\alpha )e^{-\eta (\alpha )(x-\tau )}\int_{[\lambda -x,\infty
)}e^{-z\eta (\alpha )}U_{\alpha }(dz),\text{ \ \ \ }(12)
\end{eqnarray*}%
where the first equation follows since for $x<\lambda $, $\overset{\ast }{T}%
=W_{\lambda }+W_{\tau }^{\ast }$ and upon conditioning on the sigma algebra
generated by $(W_{\lambda },I_{W_{\lambda }})$, the second equation follows
from Equation (5)\ above, the fifth equation follows from Equation (8) of
reference [7] and the last equation follows from Equation (6)\ above.
However, for $x\geq \lambda $, we have

\begin{eqnarray*}
E_{x}(e^{-\alpha \overset{\ast }{T}_{0}}) &=&E(e^{-\alpha W_{\tau }^{\ast }})
\\
&=&e^{-\eta (\alpha )(x-\tau ).}
\end{eqnarray*}%
Furthermore,

\begin{eqnarray*}
E_{\tau }(e^{-\alpha W}) &=&E_{\tau }(e^{-\alpha \overset{\ast }{T}_{0}}) \\
&=&M\eta (\alpha )\int_{[\lambda -\tau ,\infty )}e^{-x\eta (\alpha
)}U_{\alpha }(dx)\text{,}
\end{eqnarray*}%
where the last equation follows from Equation (12) above.

\ \ $\;\;\;\;\;\;\;\;\;\;\;\;\;\;\;\;\;\;\;\;\;\;\;\;\;\;\;\;\;\;\;\;\;$

Now we turn our attention to computing the cost functional for the long- run
cost average case. It follows by a Tauberian theorem that the long run
average cost per unit of time, denoted by $C(\lambda ,\tau )\;$is given by 
\[
C(\lambda ,\tau )=\frac{M[K-RE_{\tau }(E_{I_{W_{\lambda }}}(W_{\tau }^{\ast
}))]+E_{\tau }[C_{g^{\ast }}(M,I_{W_{\lambda }},\tau )]+\ C_{g}(0,\lambda
,\tau )}{E_{\tau }(W)}.\ \ (13) 
\]

\bigskip

Note that$\;E_{\tau }(W)=E_{\tau }W_{\lambda }+E_{\tau }E_{I_{W_{\lambda
}}}(W_{\tau }^{\ast })=\infty $,$\;$if \ $M\mu \leq 1$. However, if $M\mu >1$%
, we have \ \ 

\begin{eqnarray*}
E_{\tau }(W) &=&E_{\tau }(W_{\lambda })+E_{\tau }E_{I_{W_{\lambda
}}}(W_{\tau }^{\ast }) \\
&=&E_{0}(W_{\lambda -\tau })+\frac{\mu E_{\tau }(I_{_{W_{\lambda }}}-\tau
)\mu }{(\mu M-1)} \\
&=&E_{0}(W_{\lambda -\tau })+\frac{\mu E_{0}(I_{_{W_{\lambda -\tau }}})}{%
(\mu M-1)} \\
&=&E_{0}(W_{\lambda -\tau })+\frac{E_{0}(W_{\lambda -\tau })}{(\mu M-1)} \\
&=&\frac{\mu ME_{0}(W_{\lambda -\tau })}{(\mu M-1)}\ .\ \ \ \ \ (14)
\end{eqnarray*}

\ \ \ \ \ \ \ \ \ \ \ \ \ \ \ \ \ \ \ \ \ \ \ \ \ \ \ \ \ \ \ \ \ \ \ \ \ \
\ \ \ \ \ \ \ \ \ \ \ \ \ \ \ \ \ \ \ \ \ \ \ \ \ \ \ \ \ \ \ \ \ \ \ \ \ \
\ \ \ \ \ \ \ \ \ \ \ \ \ \ \ \ \ \ \ \ \ \ \ \ \ \ \ \ \ \ \ \ \ \ \ \ \ \
\ \ \ \ \ \ \ \ \ \ \ \ \ \ \ \ \ \ \ \ \ \ \ \ \ \ \ \ \ \ \ \ \ \ \ \ \ \
\ \ \ \ \ \ \ \ \ \ \ \ \ \ \ \ \ \ \ \ \ \ \ \ \ \ \ \ \ \ \ \ \ \ \ \ \ \
\ \ \ \ \ \ \ \ \ \ \ \ \ \ \ \ \ \ \ \ \ \ \ \ \ \ \ \ \ \ \ \ \ \ \ \ \ \
\ \ \ \ \ \ \ where the second equation follows from Equation (7)\ above,
the fourth equation follows from Equation (10)\ above and $E_{0}(W_{\lambda
-\tau })$\ is given in Equation (4) above.

\bigskip

Now we turn our attention to finding the stationary distribution of the
content process. It follows that, when $M\mu >1$, the content process is
ergodic. We let $Z=\underset{t\longrightarrow \infty }{\lim }\ Z_{t}$ and $%
F(z)\ $be\ the\ distribution\ function of the process $Z$. For simplicity of
the notation, we now denote $U_{0}$ and $U_{0}^{\ast }\ $by $U\ $and $%
U^{\ast }\ $,respectively. It follows from Equations (13) and (14) above that

\begin{eqnarray*}
F(z) &=&\frac{(M\mu -1)[C_{I_{[0,z]}}^{0}(0,\tau ,\lambda
)+E_{0}[C_{[0,z]}^{\alpha }(M,I_{W_{\lambda -\tau }}+\tau ,\tau )}{M\mu
E_{0}(W_{\lambda -\tau })} \\
&=&\frac{(M\mu -1)U_{I_{[0,\lambda \wedge z-\tau
)}}(0)+E_{0}[U_{I_{[0,z-\tau )}}^{\ast }(I_{W_{\lambda }})]}{M\mu
E_{0}(W_{\lambda -\tau })}\text{. \ }(15)
\end{eqnarray*}%
\newline

\ \ \ \ \ \ \ \ \ \ \ \ \ \ \ \ \ \ \ \ \ \ \ \ \ \ \ \ \ \ \ \ \ \ \ \ \ \
\ \ \ \ \ \ \ \ \ \ \ \ \ \ \ \ \ \ \ \ \ \ \ \ \ \ \ \ \ \ \ \ \ \ \ \ \ \
\ \ \ \ \ \ \ \ \ \ \ 

\bigskip

REFERENCES

\bigskip

[1] Abdel-Hameed M.S. and\ Nakhi Y. (2006)Optimal Control of a finite dam
with diffusion input and state dependent release rates. Comp. Math. Appl.
317-324.

[2] Abdel-Hameed M.S. (2000) Optimal control of dams using $P_{\lambda ,\tau
}^{M}\;$policies and penalty cost when the input process is a compound
Poisson Process with positive drift. J. Appl. Prob. 37,408-416.

[3] Abdel-Hameed M.S. and\ Nakhi Y. (1990) Optimal control of a finite dam
using $P_{\lambda ,\tau }^{M}\;$policies and penalty cost:Total discounted
and long run average cases. J.Appl.Prob.28, 888-898.

[4] Abdel-Hameed M.S. and\ Nakhi (1991) Optimal replacement and maintenance
of systems subject to semi-Markov damage. Stoch. Proc. .Appl. 37, 141-160\
.\ 

[5] Abdel-Hameed M.S. (1987) Inspection and maintenance of Devices subject
to deterioration. Advances in Applied\ Probability , 19 , 917-931.\ \ \ \ \
\ \ \ \ 

[6] Abdel-Hameed M.S.(1984) Life distribution properties of devices subject
to a Levy wear process. Mathematics of Operations Research., 9, 606-614.\ \
\ \ 

[7] Alili L. and A. E. Kyprianou (2005) Some remarks on the first passage of
Levy processes, the American put and pasting principles. The Annals of
Applied Probability. 15, 2062-2080.

[8] Bae J., Kim S. and Lee \ E.Y. (2003) Average cost under the \ $%
P_{\lambda ,\tau }^{M}$\ \ policy in a finite dam with compound Poisson
inputs. J. Appl. Prob. 40, 519--526.

[9] Blumenthal , R. M. and Getoor, R.K. (1968) Markov Processes and
Potential Theory. Academic Press, New York.\ \ \ \ \ \ \ \ \ \ \ \ \ \ \ \ \
\ \ \ \ \ \ \ \ \ \ \ \ \ \ 

[10] Dohi T., Kaio N. and Osaki S. (1995) Optimal Control of a finite dam
with a sample path constraint. Mathl.Comput. Modelling. 22, 45-51.\ \ \ \ 

[11] Lam Yeh and Lou Jiann Hua (1987) Optimal control of a finite dam.
J.Appl.Prob.24, 196-199.

[12] Lamperti, J. (1977) Stochastic Processes:\ A survey of the Mathematical
Theory. Springer Verlag, New York\ .\ \ \ \ \ 

[13] Zuckerman, D. (1977) Two- stage output procedure of a finite dam.
J.Appl.Prob.14 , 421-425.

\end{document}